\documentclass[11pt]{article}
\usepackage[latin1]{inputenc}
\usepackage[english]{babel}
\usepackage{amssymb,amsfonts, amsmath, color}
\usepackage[margin=2.7cm]{geometry}
\usepackage{xcolor}

\newtheorem{theorem}{Theorem}[section]

\newtheorem{corollary}[theorem]{Corollary}
\newtheorem{lemma}[theorem]{Lemma}
\newtheorem{remark}[theorem]{Remark}
\newtheorem{definition}[theorem]{Definition}

\newcommand{\bcl}{\begin{center}}
\newcommand{\ecl}{\end{center}}
\newcommand{\brl}{\begin{right}}
\newcommand{\erl}{\end{right}}
\newcommand{\ben}{\begin{enumerate}}
\newcommand{\een}{\end{enumerate}}
\newcommand{\overliner}{\begin{array}}
\newcommand{\earr}{\end{array}}
\newcommand{\btab}{\begin{tabular}}
\newcommand{\etab}{\end{tabular}}
\newcommand{\bdoc}{\begin{document}}
\newcommand{\edoc}{\end{document}}
\newcommand{\beqy}{\begin{eqnarray}}
\newcommand{\eeqy}{\end{eqnarray}}

\newcommand{\beqi}{\begin{eqnarray*}}
\newcommand{\eeqi}{\end{eqnarray*}}
\newcommand{\bitem}{\begin{itemize}}

\newcommand{\eitem}{\end{itemize}}
\newcommand{\nln}{\newline}
\newcommand{\newt}{\newtheorem}

\newcommand{\pa}{\partial}
\newcommand{\re}{{I\!\!R}}
\newcommand{\ren}{\re^N}
\newcommand{\xr}{x\in\re }
\newcommand{\x}{\times}
\newcommand{\dyle}{\displaystyle}
\newcommand{\ene}{{I\!\!N}}
\newcommand{\irn}{\int\limits_{\re^N}}
\newcommand{\io}{\int\limits_{\O}}
\newcommand{\meas}{{\rm meas\,}}

\newcommand{\sign}{{\rm sign}}
\newcommand{\map}{\longrightarrow }
\newcommand{\imp}{\Longrightarrow }
\renewcommand{\div}{\nabla\cdot }
\newcommand{\sen}{{\rm sen\,}}
\newcommand{\tg}{{\rm tg\,}}
\newcommand{\arcsen}{{\rm arcsen\,}}
\newcommand{\arctg}{{\rm arctg\,}}
\newcommand{\supp}{{\textsl supp\ }}
\newcommand{\ity}{\int_{-\iy}^{+\iy}}
\newcommand{\limit}{\lim\limits}
\newcommand{\limi}{\limit_{n\to\infty}}
\newcommand{\sumi}{\sum\limits_{n=1}^{\infty}}
\newcommand{\ulu}{\underline u}
\newcommand{\ulw}{\underline w}
\newcommand{\ulz}{\underline z}
\newcommand{\ulv}{\underline v}
\newcommand{\uls}{\underline s}
\newcommand{\olu}{\overline u}
\newcommand{\olv}{\overline v}
\newcommand{\ols}{\overline s}
\newcommand{\ob}{\overline\b}
\newcommand{\ovar}{\overline\var}
\newcommand{\wv}{\widetilde v}
\newcommand{\wu}{\widetilde u}
\newcommand{\ws}{\widetilde s}
\renewcommand{\a }{\alpha }
\renewcommand{\b }{\beta }
\newcommand{\g }{\gamma}
\newcommand{\G }{\Gamma }
\renewcommand{\d }{\delta }

\newcommand{\D }{\Delta }
\newcommand{\e }{\varepsilon }
\newcommand{\z }{\zeta }
\renewcommand{\l }{\lambda }
\renewcommand{\L }{\Lambda }
\newcommand{\m }{\mu }
\newcommand{\n }{\nabla }
\newcommand{\s }{\sigma }
\newcommand{\Sig }{\Sigma }
\renewcommand{\t }{\tau }
\newcommand{\var }{\varphi }
\renewcommand{\o }{\omega }
\renewcommand{\O }{\Omega }
\newcommand{\bR}{{\bf R}}
\newcommand{\bC}{{\bf C}}
\newcommand{\bZ}{{\bf Z}}
\newcommand{\bN}{{\bf N}}
\newcommand{\bQ}{{\bf Q}}
\newcommand{\bK}{{\bf K}}
\newcommand{\bI}{{\bf I}}
\newcommand{\bv}{{\bf v}}
\newcommand{\bV}{{\bf V}}
\DeclareMathOperator{\suppo}{supp}

\def\qed{\unskip\kern 6pt \penalty 500
\raise -2pt\hbox{\vrule \vbox to10pt{\hrule width 4pt
\vfill\hrule}\vrule}\par}

\newenvironment{Proof}{\removelastskip\vskip12pt
plus 1pt \noindent\em\rm}{\hfill {\qed \hskip .2cm}}

\title
{Conditions at infinity for the inhomogeneous filtration equation}
\author{Gabriele Grillo\thanks{Dipartimento di Matematica, Politecnico di Milano, Piazza Leonardo da Vinci 32,
20133 Milano, Italy (gabriele.grillo@polimi.it).}\,\,, Matteo
Muratori\thanks{Dipartimento di Matematica, Politecnico di Milano,
Piazza Leonardo da Vinci 32, 20133 Milano, Italy
(matteo1.muratori@mail.polimi.it).}\,\, and Fabio
Punzo\thanks{Dipartimento di Matematica "G. Castelnuovo",
Universit\`a di Roma ``La Sapienza'', Piazzale Aldo Moro 5, 00185
Roma, Italy (punzo@mat.uniroma1.it).}}

\date{}

\begin{document}
\maketitle

\abstract{We investigate existence and uniqueness of solutions to the filtration equation with an inhomogeneous density in ${\mathbb R}^N$ ($N\ge3$), approaching at infinity a given continuous datum of Dirichlet type.}


\bigskip
\smallskip

\section{Introduction}\setcounter{equation}{0}
We provide sufficient conditions for existence and uniqueness of
{\it bounded} solutions to the following
nonlinear Cauchy problem (given $T>0$):
\begin{equation}\label{e1}
\left\{
\begin{array}{ll}
\,   \rho\, \pa_t u = \Delta\big[G(u)\big]
&\textrm{in}\,\,\re^N\times (0, T]=:S_T
\\& \\
\textrm{ }u \, = u_0& \textrm{in\ \ } \re^N\times \{0\} \,,
\end{array}
\right.
\end{equation}
where $\rho=\rho(x)$ does not depend on $t$. Concerning the density $\rho$, the initial
condition $u_0$ and the nonlinearity $G$ we shall mostly assume the following:
\[
\textrm{\ \ } \left\{
\begin{array}{l}
(i) \quad \;  \rho\in C(\re^N), \, \rho>0 \, ;
\\
(ii)\quad \! G \in C^1(\re) , \, G(0)=0,\, G'(s)>0\,\hbox{\,for any}\,\,s\in\re\setminus\{0\}, \\
\qquad \,\, G'\,\hbox{decreasing in}\, (-\d,0)\,\hbox{and
increasing in}\,\,(0,\d) \\
\qquad \,\,\hbox{if}\,\, G'(0)=0\,\,(\d>0), \, \hbox{for some $\delta>0$} ;
\\ (iii)\;\; u_0\in L^{\infty}(\re^N)\cap C(\re^N) .
\end{array}
\right. \leqno(H_0) \]
A typical choice for the function $G$ is
$G(u)=|u|^{m-1}u$ for some $m\ge 1.$ In this case, for $m>1$, the
differential equation in problem \eqref{e1} becomes the {\it
inhomogeneous porous media equation}, which arises in various
situations of physical interest. We quote, without any claim of
generality, the papers \cite{KR1}, \cite{KR}, \cite{Eid90},
\cite{RV}, \cite{DG+08}, \cite{DNS}, \cite{RV08}, \cite{RV09}, \cite{P1}--\cite{P2},
\cite{KRV10}, \cite{GM}, \cite{GMP} as references
on this topic, and the recent monograph \cite{Vaz07} as a general
reference on the porous media equation.

As it is well-known, if assumption $(H_0)$ is satisfied, then there
exists a bounded solution of problem \eqref{e1} (see, $e.g.$, \cite{KR}, \cite{EK}, \cite{RV}). Moreover, if $N=1$ or $N=2$, and
$\rho\in L^\infty(\re^N)$, then the solution to problem \eqref{e1}
is unique (see \cite{GHP}).

When $N\ge 3$, we can have uniqueness or nonuniqueness of bounded
solutions to problem \eqref{e1}, in dependence of the behavior
{\it at infinity} of the density $\rho$. In fact, given $R>0$, set $B_R:=\{x\in
\re^N: \, |x|< R \}$ and $B_R^c:= \re^N \setminus B_R$. Suppose that $\rho$ does not decay too fast at infinity, in the sense that
there exist $\widehat{R}>0$ and $\underline{\rho} \in C([\widehat{R}, \infty))$ such that
$\rho(x)\ge \underline{\rho}(|x|)>0$ for
all $ x \in B_{\widehat{R}}^c$, with $\int_{\widehat{R}}^\infty \eta \underline{\rho}(\eta)\,
d\eta\,=\, \infty$. Then problem \eqref{e1} admits \emph{at most} one
bounded solution (see \cite{P1}, \cite{P2}). A natural choice for $\underline{\rho}$ above is $\underline{\rho}(\eta):=\eta^{-\a}$ ($\eta\in
[\widehat{R},\infty)$) for some $\a\in (-\infty, 2]\,$ and
$\widehat{R}>0$.

On the contrary if $\rho $ decays sufficiently fast at infinity, in the sense that $\Gamma\ast\rho \in L^{\infty}(\re^N)$, where $\Gamma$ is the fundamental solution of the Laplace equation
in $\re^N$, then nonuniqueness prevails (see \cite{P1}, \cite{P2} and also \cite{KPT} for the linear case, namely $G(u)=u$). To be specific, for
any function $A\in \textrm{Lip}\,[0, T]$ with $A(0)=0$ there exists a
solution $u$ of problem \eqref{e1} such that
\begin{equation}\label{e5}
\lim_{R\to \infty} \frac{1}{|\pa B_R|}\int_{\pa B_R}\big| U(x, t)
- A(t)\big|\, d\s\,=\,0\,
\end{equation}
uniformly with respect to $t\in [0,T]$, where
\begin{equation}\label{e3} U(x, t):= \int_0^t G(u(x,\t))\,
d\t  \ \ \ \forall(x,t) \in S_T \,.
\end{equation}
The condition $\Gamma\ast\rho \in L^{\infty}(\re^N)$ can be replaced by the following stronger (but more explicit) condition: there exists $\widehat{R}>0$ and $\overline\rho\in C([\widehat{R}, \infty))$ such that $\rho(x)\le \overline{\rho}(|x|)$ for
all $x\in B_{\widehat{R}}^c$, with $\int_{\widehat{R}}^\infty \eta \overline{\rho}(\eta)\,
d\eta\,<\, \infty \,$. Then, instead of \eqref{e5}, we can impose
that
\begin{equation}\label{e5a}
\lim_{|x|\to \infty}\big| U(x, t) - A(t) \big|\,=\, 0\,
\end{equation}
uniformly with respect to $t\in [0,T]$ (which clearly implies \eqref{e5}), with $U$ defined in
\eqref{e3}. A natural choice for $\overline\rho$ as above is $\overline{\rho}(\eta):=\eta^{-\a}$ ($\eta\in
[\widehat{R},\infty)$) for some $\a\in (2, \infty]\,$ and
$\widehat{R}>0$.

\smallskip

Observe that equalities \eqref{e5} and \eqref{e5a} can be also
regarded as nonhomogeneous Dirichlet conditions at infinity  in a suitable integral sense. From this point of view, it seems natural to study whether imposing conditions at infinity \emph{in a pointwise sense} that resembles more closely the usual Dirichlet boundary conditions restores existence and uniqueness of solutions. In fact, up to now, it was only known that there exists at most one solution $u\in L^\infty(S_T)$
to problem \eqref{e1} satisfying condition \eqref{e5} or
\eqref{e5a} either when $G(u)=u$ (see the important results obtained, in such linear case, in \cite{KPT}) or when $u_0\ge
0$ and $A\equiv 0$ (see \cite{EK}). However, the methods used
to obtain the mentioned uniqueness results did not work for general $G$ and $A$.

In this paper we shall then address existence and uniqueness of bounded
solutions to problem \eqref{e1} satisfying at infinity suitable nonhomogeneous Dirichlet conditions in a pointwise sense. More precisely, at first we shall prove that if $\rho$ decays sufficiently fast at infinity, the diffusion is non-degenerate in an appropriate sense, $u_0\in C(\re^N)$ and $\lim_{|x|\to \infty} u_0(x)$ exists and is finite then for any $a\in C([0,T])$
with
\begin{equation}\label{e65}
a(0)=\lim_{|x|\to \infty} u_0(x)
\end{equation}
there exists a bounded solution $u$ to problem \eqref{e1} satisfying
\begin{equation}\label{e63}
\lim_{|x|\to \infty} u(x,t)=a(t)\quad \textrm{uniformly for}\,\,\,
t\in [0,T]
\end{equation}
(see Theorem \ref{t2}). Furthermore, we can remove the assumption of nondegeneracy of the diffusion
for suitable classes of initial data $u_0$
and conditions at infinity $a$. Indeed if $a_0:=\lim_{|x|\to
\infty} u_0(x)$ exists and is finite and $(H_0)$ holds true, then there
exists a bounded solution to problem \eqref{e1} such that
\begin{equation}\label{e64}
\lim_{|x|\to \infty} u(x,t)=a_0 \quad \textrm{uniformly for}\,\,\,
t\in [0,T]\,
\end{equation}
(see Theorem \ref{t3} and Remark \ref{osst3}).

Moreover, if $(H_0)$ holds true and $\rho$ decays sufficiently fast at infinity, then there exists a bounded solution $u$ to problem \eqref{e1} satisfying \eqref{e63} for any
$a\in C([0,T])$ with $a>0$ in $[0,T]$, provided $u_0$ complies with
\eqref{e65} (see Theorem \ref{t4}).

Let us explain that  in \cite{P1}, generalizing arguments used in \cite{KPT},  the prescription of
conditions \eqref{e5}  for solutions to \eqref{e1} is made by constructing
suitable {\it barriers} at infinity, that are sub-- or
supersolutions to appropriate associated linear elliptic problems.
Instead, in the present case to impose at infinity Dirichlet
conditions in a pointwise sense we will construct, in a
neighborhood of each $t_0\geq 0$, suitable time-dependent {\it
barriers} at infinity, that are sub-- or supersolutions to proper
associated nonlinear parabolic problems.

Actually, in the existence results, hypothesis \eqref{e65} can be removed, upon requiring that the Dirichlet condition at infinity is attained uniformly for $t\in [\t, T]$, for any $0<\t<T$ (see Remark \ref{texib})  and, in the degenerate case, a further technical condition holds.

\smallskip

Finally, we shall prove that the weaker condition
\begin{equation}\label{eq: e63Weak}
\lim_{|x|\to \infty} u(x,t)=a(t) \ \ \ \forall
t\in [0,T]
\end{equation}
implies uniqueness for general $G$ satisfying $(H_0)(ii)$, bounded $\rho$ and $a\in
C([0,T])$ (see Theorem \ref{t1}). Arguments used in proving
uniqueness are modeled after those in \cite{ACP} (where
$\rho\equiv 1, N=1$) and \cite{GHP} (where $N=2$), for cases in
which uniqueness was proved in the class of solutions not
satisfying additional conditions at infinity. Although this is not
our case, we use an analogous strategy, combined with the fact
that solutions attain a datum at infinity in a pointwise sense.
This permits to conclude.

\noindent We thank the referees for their careful reading of the original version of this manuscript. In particular we thank one of them for pointing out that our arguments could be modified to yield the conclusions discussed in Remark \ref{texib}, the other one for some
suggestions that have improved the presentation.

\section{Existence and uniqueness results}\setcounter{equation}{0}
Solutions, sub- and supersolutions to problem \eqref{e1} are always meant in the following sense.

\begin{definition}\label{defsol}
By a {\em solution} to problem \eqref{e1} we mean a function $u\in
C(S_T)\cap L^\infty(S_T)$ such that
\begin{equation}\label{e2}
\begin{split}
\int_0^\t\int_{\Omega_1}\big\{\rho\, u\, \pa_t \psi + G(u) \Delta
\psi \big\}\, d x d t = & \int_{\Omega_1}
\rho\,\big[u(x,\t)\psi(x,\t) - u_0(x) \psi(x, 0)\big]\,dx \\
& + \int_0^\t\int_{\pa \Omega_1} G(u) \langle\nabla \psi,\, \nu\rangle
\,d\s dt\hspace{2 cm}
\end{split}
\end{equation}
for any bounded open set\, $\Omega_1\subseteq \re^N$ with smooth
boundary $\pa\Omega_1,\,\t\in (0,T],\,\psi\in
C^{2,1}(\overline{\Omega}_1 \times[0,\t]), \psi\ge 0, \psi=0$ in
$\pa \Omega_1\times [0,\t]$. Here $\nu$ denotes the outer normal
to $\Omega_1$ and $\langle\cdot, \cdot \rangle$ the scalar product
in $\re^N$.

{\em Supersolutions (subsolutions)} to \eqref{e1} are defined
replacing $``="$ by $``\le"$\, ($``\ge"$, respectively) in
\eqref{e2}.
\end{definition}

These kind of solutions are sometimes referred to as \it very weak
\rm solutions. Observe that, according to Definition \ref{defsol},
solutions to problem \eqref{e1} we deal with are in
$S_T$.

\subsection{Existence}
In the case of \emph{nondegenerate} nonlinearities, we have the following
result.\normalcolor
\begin{theorem}\label{t2}
Let $N\ge 3$. Assume that $\rho\in C(\re^N)$, $\rho>0$, $G\in C^1(\re)$ with
$G(0)=0$, $G'(s)\ge\a_0>0$ for any $s\in\re$ and $u_0\in C(\re^N)$ with $\lim_{|x|\to \infty}
u_0(x)$ existing and being finite. Assume also that there exist $\widehat{R}>0$ and
$\overline\rho\in C([\widehat{R}, \infty))$ such that $\rho(x)\le \overline{\rho}(|x|)$ for
any $x\in B_{\widehat{R}}^c$, with $\int_{\widehat{R}}^\infty \eta \overline{\rho}(\eta)\,
d\eta\,<\, \infty$.

Finally, let $a\in C([0,T])$ and suppose that
\[a(0)=\lim_{|x|\to \infty} u_0(x) \, . \]
Then there exists a solution to problem \eqref{e1} such that
\[ \lim_{|x|\to \infty} u(x,t)=a(t) \quad \textnormal{uniformly for}\
t\in [0,T] \, . \]
\end{theorem}

For appropriate classes of data and possibly \emph{degenerate} nonlinearities of porous media type, we shall prove the following results.

\begin{theorem}\label{t3}
Let $N\ge 3$. Let assumption $(H_0)$ be satisfied. Suppose that
\begin{equation}\label{e66}
\lim_{|x|\to \infty} u_0(x)\,=\, a_0
\end{equation}
for some $a_0\in\re.$ Then there exists a solution to problem
\eqref{e1} such that
\[ \lim_{|x|\to \infty} u(x,t)=a_0 \quad \textnormal{uniformly for}\
t\in [0,T] \, . \]
\end{theorem}

\begin{remark}\label{osst3} Let assumption $(H_0)$ be satisfied and suppose that $\rho$ does not decay too fast at infinity in the sense  that there exist $\widehat{R}>0$ and
$\underline\rho\in C([\widehat{R}, \infty))$ such that $\rho(x)\ge \underline{\rho}(|x|)>0$ for
any $x\in B_{\widehat{R}}^c$, with $\int_{\widehat{R}}^\infty \eta \underline{\rho}(\eta)\,
d\eta\,=\infty$. Assume also that \eqref{e66} holds. Then by the uniqueness result recalled in the Introduction, and by Theorem \ref{t3}, the \emph{unique} solution to problem \eqref{e1} necessarily satisfies
\[ \lim_{|x|\to \infty} u(x,t)=a_0 \quad \textnormal{uniformly for}\
t\in [0,T] \, . \]
\end{remark}

\begin{theorem}\label{t4} Let $N\ge 3$. Let assumption $(H_0)$ be satisfied. Suppose that there exist $\widehat{R}>0$ and $\overline\rho\in C([\widehat{R}, \infty))$ such that $\rho(x)\le \overline{\rho}(|x|)$ for any $x\in B_{\widehat{R}}^c$, with $\int_{\widehat{R}}^\infty \eta \overline{\rho}(\eta)\,d\eta\,<\, \infty$. Let $a\in C([0,T])$, with $a(t)>0$ for all $t \in [0,T]$. Assume also that
\[a(0)=\lim_{|x|\to \infty} u_0(x) \, .\]
\noindent Then there exists a solution to problem \eqref{e1} such that
\[ \lim_{|x|\to \infty} u(x,t)=a(t) \quad \textnormal{uniformly for}\
t\in [0,T] \, . \]
\end{theorem}

\begin{remark}\label{texib}
\begin{itemize}
\item[$(i)$] In Theorem \ref{t2}, if we \emph{do not} assume that $a(0)=\lim_{|x|\to \infty} u_0(x)$, then the conclusion remains true, replacing the property $\lim_{|x|\to \infty} u(x,t)=a(t)$ uniformly for any $t\in [0,T]$ by the following:
\begin{equation}\label{e63b}
\textrm{for any}\;\,\t\in (0,T) \, ,\;\;   \lim_{|x|\to \infty} u(x,t) =
a(t) \quad \textnormal{uniformly for} \;\, t\in [\tau, T]\,.
\end{equation}
\item[$(ii)$] In Theorem \ref{t4}, if we \emph{do not} assume that $a(0)=\lim_{|x|\to \infty} u_0(x)$, then the conclusion remains true, provided we replace the property $\lim_{|x|\to \infty} u(x,t)=a(t)$ uniformly for any $t\in [0,T]$ by \eqref{e63b} and we also require that
\begin{equation}\label{er4}
I:=\inf_{x \in B_{R_0}^c}u_0(x) \, , \ S:=\sup_{t \in (0,\epsilon) } a(t)  \, ,  \   2 \, G (I) > G\left( S \right)
\end{equation}
for some $R_0,\epsilon >0$. Clearly, \eqref{er4} is technical and is needed to make our proof hold under more general assumptions. We do not know whether the result is still valid without assuming it, but notice that \eqref{er4} certainly holds if $I$ is large enough compared to $S$, so that possible problems occur only if the initial datum is, in a suitable sense, small at infinity.
\end{itemize}
\end{remark}

See the end of Section \ref{secpexi} for comments on the minor changes needed in the proof of the corresponding theorems to obtain statements $(i)-(ii)$.

\begin{remark}\label{osst4}
Note that the hypotheses made in Theorem \ref{t4} allow to assume
as initial data functions $u_0$ which may be nonpositive in some
compact subset $K\subset \re^N$.
\end{remark}
\subsection{Uniqueness}\setcounter{equation}{0}
We shall prove the following uniqueness result in the general case
of possibly degenerate nonlinearities.

\begin{theorem}\label{t1}
Let $N\ge 3$. Let assumption $(H_0)$ be satisfied, and suppose that $a\in L^\infty(0,T), \ \rho$ $ \in L^\infty(\re^N)$. Then there exists \emph{at most} one solution to problem \eqref{e1} such that
$$ \lim_{|x|\to \infty} u(x,t)=a(t)\ \ \ \textnormal{for almost every} \ t \in (0,T) \, . $$
\end{theorem}
From Theorems \ref{t3} and \ref{t1} we deduce the following.
\begin{corollary}\label{cor2}
Let $N\ge 3.$ Let assumption $(H_0)$ be satisfied, and suppose
that $\rho \in L^\infty(\re^N)$. Then there exists a unique
solution to problem \eqref{e1} such that
\[
\lim_{|x|\to \infty} u(x,t)=a_0 \quad \textnormal{uniformly for}\,\,
t\in [0,T].
\]
\end{corollary}

\begin{remark}\label{ossc2}
When $(H_0)$ is satisfied, $\rho$ belongs to $L^\infty(\re^N)$ and fulfils the assumptions appearing in Remark \ref{osst3}, then the conclusion of Corollary \ref{cor2} is in agreement with such Remark.
\end{remark}
As a consequence of Theorems \ref{t4} and \ref{t1} we get
\begin{corollary}\label{cor3}
Let $N \ge 3$. Let the assumptions of Theorem \ref{t4} be satisfied, and suppose that $\rho \in L^\infty(\re^N)$. Then there exists a unique solution to problem \eqref{e1} such that
\[
\lim_{|x|\to \infty} u(x,t)\,=\, a(t) \quad \textnormal{uniformly
for}\,\, t\in [0,T] \, .
\]
\end{corollary}

Finally, in the case of \it nondegenerate \rm nonlinearities, from
Theorems \ref{t2} and \ref{t1} we obtain the following.
\begin{corollary}\label{cor1}
Let $N \ge 3$. Let the assumptions of Theorem \ref{t2} be satisfied, and suppose that $\rho \in L^\infty(\re^N)$. Then there exists a unique solution to problem \eqref{e1} such that
\[
\lim_{|x|\to \infty} u(x,t)\,=\, a(t) \quad  \textnormal{uniformly
for}\,\, t\in [0,T] \, .
\]
\end{corollary}

\section{Existence: proofs}\label{secpexi}
\setcounter{equation}{0}
In view of  the assumptions on $\rho$ made in the hypotheses of Theorem \ref{t2} or \ref{t4}, there exists a function $V=V(|x|)\in C^2(B_{\widehat R}^c)$ such that
\begin{equation}\label{e73}
\Delta V\, \le\, - \rho \quad \textrm{in}\;\; B_{\widehat R}^c \,,
\end{equation}
\[
V(|x|)\,>\,0 \quad  \forall x\in B_{\widehat
R}^c \, ,
\]
\[|x|\mapsto V(|x|)\;\; \textrm{is nonincreasing} \,,\]
\begin{equation}\label{e71}
\lim_{|x|\to \infty} V(x)\,=\,0\,;
\end{equation}
here $\widehat{R}>0$ can be assumed to be equal to the one that appears in the hypotheses of Theorem \ref{t2} or \ref{t4}. 

In some of the forthcoming proofs we shall make use of the function $G^{-1}$, whose domain $D$ need not coincide with ${\mathbb R}$. As we are dealing with bounded data $u_0$ (and, by the maximum principle, with bounded solutions), this makes no problem since one can modify the definition of $G(x)$ for $|x|$ large so that such a function is a bijection from ${\mathbb R}$ to itself, without changing the evolution of $u_0$.

\smallskip

Hereafter, for any $j\in \ene$, $\zeta_j$ will always be a function having the following properties: $\zeta_j\in
C^{\infty}_c(B_j)$ with $0\le \zeta_j\le 1$ and $ \zeta_j\equiv 1$ in $B_{j/2}$.

\bigskip

\noindent{\it Proof of Theorem \ref{t2}\,.} Since $a\in C([0,T])$
and $G\in C^1(\re)$ is increasing, for any $t_0\in [0,T]$, $\s>0$ there exists $\d=\d(\s)>0$ such that
\begin{equation}\label{e72}
G^{-1}\big[G(a(t_0))-\s\big]\, \le\, a(t)\,\le
G^{-1}\big[G(a(t_0))+\s\big]\quad \forall t\in
[\underline t_\d, \overline t_\d]\,,
\end{equation}
where $\underline t_\d:= \max\{t_0-\d, 0\}$ and $\overline t_\d:=
\min\{ t_0+\d, T\}$. Moreover, in view of the assumptions on $u_0$, for
any $\s>0$ there exists $R=R(\s)>\widehat R$ such that
\begin{equation}\label{e72a}
G^{-1}\big[G(a(0))-\s\big]\, \le\, u_0(x)\,\le
G^{-1}\big[G(a(0))+\s\big]\quad \forall x\in
B_R^c \, .
\end{equation}
For any $j\in \ene$, let $u_j\in C\big(\overline{B}_j\times
[0,T]\big)$ be the unique solution (see, $e.g.$, \cite{LSU}) to the
problem
\begin{equation}\label{e68}
\left\{
\begin{array}{ll}
\,   \rho\, \pa_t u = \Delta\big[G(u)\big] &\textrm{in}\,\,B_j\times (0, T) \\&\\
\, u\,=\, a(t) &\textrm{on}\,\, \pa B_j\times (0,T)
\\ & \\ \, u \, = u_{0,j} &\textrm{in}\,\, \overline{B}_j\times \{0\}\,,
\end{array}
\right.
\end{equation}
where
\[u_{0,j}:= \zeta_j \, u_0 + (1-\zeta_j) \, a(0) \quad \textrm{in} \;\; \overline{B}_j\,.\]
By comparison principles,
\begin{equation}\label{e69}
|u_j|\le K:=\max\{\|u_0\|_\infty, \|a\|_\infty\}\quad
\textrm{in}\;\;B_j\times (0,T)\,.
\end{equation}
It is a matter of usual compactness arguments (see, $e.g.$, \cite{LSU}) to show that there
exists a subsequence $\{u_{j_k}\}\subseteq\{u_j\}$ which converges, as $k\to \infty$, \emph{locally uniformly} in $\re^N\times (0,T)$ to a solution $u$ to problem \eqref{e1}.

Hence, it remains to prove that
\[
\lim_{|x|\to \infty} u(x,t)=a(t)  \quad  \textnormal{uniformly
for}\,\, t\in [0,T] \, .
\]
To this end, let $t_0\in [0,T]$. Define
\[\underline{w}(x,t):=G^{-1}\big[-\underline{M} V(x) -\s +G(a(t_0)) -\underline{\l} (t-t_0)^2\big] \quad \forall(x,t) \in  B_{\widehat R}^c \times (\underline t_\d, \overline t_\d) \, , \]
where $\underline{M}>0$ and $\underline{\l}>0$ are constants to be chosen later. By
the assumptions and \eqref{e73},
\begin{equation}\label{e74}
\rho(x)\pa_t \underline w-\Delta[G(\underline w)] =  -\rho(x)\frac{2\underline{\l}
(t-t_0)}{G'(w)} + \underline{M} \Delta V   \le  \rho(x)\left(\frac{2\l
\d}{\a_0} - \underline{M} \right)\le 0
\end{equation}
$$ \textrm{in} \; \; B_{\widehat R}^c \times (\underline t_\d, \overline t_\d)  $$
providing that
\begin{equation}\label{e74a}
\underline{M} \ge \frac{2 \underline{\l} \d}{\a_0}\,.
\end{equation}
For any $j\in \ene$, $j>R$, let
$$N_{R,j}:= B_j\setminus \overline{B}_R \, ,$$
$R$ being as in \eqref{e72a}. We have
\begin{equation}\label{e75}
\underline w(x,t)\le -K \quad \forall(x,t) \in \pa
B_R\times (\underline t_\d , \overline t_\d)
\end{equation}
provided
\begin{equation}\label{e75a}
\underline{M} \ge \frac{G(\|a\|_\infty)- G(-K)}{V(R)}\,.
\end{equation}
Furthermore,
\begin{equation}\label{e76}
\underline w(x,t)\le G^{-1}\big[G(a(t_0)) -\s\big]\normalcolor \quad
\forall(x,t)\in \pa B_j\times (\underline t_\d,
\overline t_\d).
\end{equation}
When $\underline t_\d=0$ there holds
\begin{equation}\label{e76b}
\underline w(x,t)\le G^{-1}\big[G(a(t_0)) -\s\big]\quad
\forall (x,t) \in \overline{N}_{R,j}\times \{\underline t_\d\}\,,
\end{equation}
whereas when $\underline t_\d>0$ we have
\begin{equation}\label{e78}
\underline w(x,t)\le G^{-1}\big[G(a(t_0)) -\underline{\l} \d^2\big]\le - K
\quad \forall (x,t)\in  \overline{N}_{R,j}\times \{\underline
t_\d\}
\end{equation}
provided
\begin{equation}\label{e78a}
\underline{\l} \ge \frac{G(\|a\|_\infty)- G(-K)}{\d^2}\,.
\end{equation}
Suppose that conditions \eqref{e74a}, \eqref{e75a} and
\eqref{e78a} are satisfied. Hence, from \eqref{e74} and \eqref{e75}--\eqref{e78} we infer that $\underline w$ is a
subsolution to the problem
\begin{equation}\label{e79}
\left\{
\begin{array}{ll}
\,  \rho\, \pa_t u = \Delta\big[G(u)\big] &\textrm{in}\,\,N_{R,j}\times (\underline t_\d, \overline t_\d) \\&\\
u\, =\, -K &\textrm{on}\,\,\pa B_{R}\times (\underline t_\d, \overline t_\d)\\&\\
u\,=\,G^{-1}\big[G(a(t_0))-\s\big] &\textrm{on}\,\, \pa
B_j\times(\underline t_\d, \overline t_\d)
\\& \\ \,\! u \, = - K &\textrm{in}\,\, \overline{N}_{R,j}\times \{\underline t_\d\}
\end{array}
\right.
\end{equation}
when $\underline t_\d>0$, whereas it is a subsolution to the problem
\begin{equation}\label{e80}
\left\{
\begin{array}{ll}
\,   \rho\, \pa_t u = \Delta\big[G(u)\big] &\textrm{in}\,\,N_{R,j}\times (\underline t_\d, \overline t_\d) \\&\\
u\, =\, -K &\textrm{on}\,\,\pa B_{R}\times (\underline t_\d, \overline t_\d)\\&\\
u\,=\,G^{-1}\big[G(a(t_0))-\s\big] &\textrm{on}\,\, \pa
B_j\times(\underline t_\d, \overline t_\d)
\\& \\ \,\! u \, = G^{-1}\big[G(a(t_0))-\s\big] &\textrm{in}\,\, \overline{N}_{R,j}\times \{\underline t_\d\}
\end{array}
\right.
\end{equation}
when $\underline t_\d=0$.

\smallskip
On the other hand, \eqref{e72}, \eqref{e72a} and \eqref{e69} show
that the boundary data for the solutions to \eqref{e68} and \eqref{e79}, \eqref{e80} are correctly ordered on each part of the parabolic boundary of $N_{R,j} \times (\underline t_\d, \overline t_\d)$. In particular, we
 deduce that $u_j$ is a supersolution to problem \eqref{e79}
when $\underline t_\d>0$, while it is a supersolution to problem
\eqref{e80} when $\underline t_\d=0$. Therefore, by comparison principles,
\begin{equation}\label{e81}
\underline w\le u_j\quad \textrm{in}\;\; N_{R,j}\times (\underline
t_\d, \overline t_\d)\,.
\end{equation}

Now let us define
\[\overline w(x,t):=G^{-1}\big[\overline{M} V(x) +\s +G(a(t_0)) +\overline{\l} (t-t_0)^2\big] \quad \forall(x,t)  \in B_{\widehat R}^c \times (\underline t_\d, \overline t_\d) \,,\]
with
\[
\overline{M} \ge \max \left\{\frac{2\overline{\l} \d}{\a_0},
\frac{G(K)-G(-\|a\|_\infty)}{V(R)}\right\},
\]
and
\[
\overline{\l} \ge \frac{G(K)-G(-\|a\|_\infty)}{\d^2}\,.
\]
By arguments analogous to those used above, we can infer that $\overline w$ is a supersolution to the problem
\begin{equation}\label{e83}
\left\{
\begin{array}{ll}
\,   \rho\, \pa_t u = \Delta\big[G(u)\big] &\textrm{in}\,\,N_{R,j}\times (\underline t_\d, \overline t_\d) \\&\\
u\, =\, K &\textrm{on}\,\,\pa B_{R}\times (\underline t_\d, \overline t_\d)\\&\\
u\,=\,G^{-1}\big[G(a(t_0))+\s\big] &\textrm{on}\,\, \pa
B_j\times(\underline t_\d, \overline t_\d)
\\& \\ \,\! u \, =  K &\textrm{in}\,\, \overline{N}_{R,j}\times \{\underline t_\d\}
\end{array}
\right.
\end{equation}
when $\underline t_\d>0$, whereas it is a supersolution to the problem
\begin{equation}\label{e84}
\left\{
\begin{array}{ll}
\,   \rho\, \pa_t u = \Delta\big[G(u)\big] &\textrm{in}\,\,N_{R,j}\times (\underline t_\d, \overline t_\d) \\&\\
u\, =\, K &\textrm{on}\,\,\pa B_{R}\times (\underline t_\d, \overline t_\d)\\&\\
u\,=\,G^{-1}\big[G(a(t_0))+\s\big] &\textrm{on}\,\, \pa
B_j\times(\underline t_\d, \overline t_\d)
\\& \\ \,\! u \, = G^{-1}\big[G(a(t_0))+\s\big] &\textrm{in}\,\, \overline{N}_{R,j}\times \{\underline t_\d\}
\end{array}
\right.
\end{equation}
when $\underline t_\d=0$.

As before, from \eqref{e72}, \eqref{e72a} and \eqref{e69} we
deduce that $u_j$ is a subsolution to problem \eqref{e83}
when $\underline t_\d>0$, while it is a subsolution to problem
\eqref{e84} when $\underline t_\d=0$. By comparison principles,
\begin{equation}\label{e85}
u_j\le \overline w\quad \textrm{in}\;\; N_{R,j}\times (\underline
t_\d, \overline t_\d)\,.
\end{equation}
From \eqref{e81} and \eqref{e85} with $j=j_k$, sending $k\to
\infty$, we then obtain
\begin{equation}\label{e86}
\underline w\le u\le \overline w\quad \textrm{in}\;\;
B_R^c\times(\underline t_\d, \overline t_\d) \, .
\end{equation}
By \eqref{e86} and \eqref{e71} we get that for $|x|$ large enough, \emph{independently of} $ t_0 \in [0,T]$, there holds
\[G^{-1}\big[G(a(t_0)) -2\s\big] \le u(x,t_0)
\le G^{-1}\big[G(a(t_0))+2\s\big] \,. \]
In order to complete the proof one just lets $\s \to 0^+$. \hfill $\square$
\bigskip

\noindent{\it Proof of Theorem \ref{t3}\,.} As in the proof of the previous result note that, thanks to \eqref{e66}, for any $\s>0$ there exists
$R=R(\s)>0$ such that
\[
G^{-1}\big[G(a_0)-\s \big] \le u_0(x)\le
G^{-1}\big[G(a_0)+\s\big]\quad \forall  x \in
B_R^c \,.
\]

In view of assumption ($H_0$), by standard results (see, $e.g.$, \cite{ACP}), for any $j\in \ene$
there exists a unique solution $u_j$ to the problem
\[
\left\{
\begin{array}{ll}
\,   \rho\, \pa_t u = \Delta\big[G(u)\big] &\textrm{in}\,\,B_j\times (0, T) \\&\\
u\,=\, a_0 &\textrm{on}\,\, \pa B_j\times (0,T)
\\& \\ \,\! u \, = u_{0,j} &\textrm{in}\,\, \overline{B}_j\times \{0\}\,,
\end{array}
\right.
\]
where
\[u_{0,j}:= \zeta_j u_0 + (1-\zeta_j) \, a_0 \quad \textrm{in} \;\; \overline B_j\,.\]
Note that, by the results of \cite{DiB3}, $u_j\in C\big(\overline
B_j\times [0,T]\big)$. By comparison principles,
\[
|u_j|\le K:=\max\{ \|u_0\|_\infty, |a_0| \} \quad
\textrm{in}\;\;B_j\times (0,T)\,.
\]
By usual compactness techniques (one can use \cite[Lemma
5.2]{DiB1} and a diagonal argument), there exists a subsequence
$\{u_{j_k}\}\subseteq\{u_j\}$ which converges, as $k\to \infty,$
locally uniformly in $\re^N\times (0,T)$ to a solution $u$ to problem \eqref{e1}.

\smallskip
Let
\[\Gamma(x)\equiv \Gamma(|x|):=|x|^{2-N} \quad \forall x \in \re^N\setminus \{0\} \,.\]
Clearly,
\[
\Delta \Gamma\,=\,0\quad \textrm{in}\;\; \re^N\setminus \{0\}\,,
\]
\[
\Gamma\,>\,0\quad \textrm{in}\;\; \re^N\setminus \{0\} \, ,
\]
\begin{equation}\label{e89}
\lim_{|x|\to \infty}\Gamma(|x|)\,=\, 0\,.
\end{equation}

\noindent Define
\[\overline W(x):=G^{-1}\big[\overline{M} \, \Gamma(x) + \s + G(a_0) \big] \quad \forall x\in \re^N\setminus\{0\}\,,
\]
where
\begin{equation}\label{e95}
\overline{M} \ge \frac{G(K)-G(a_0)}{\Gamma(R)} \, .
\end{equation}
Then
\begin{equation}\label{e93}
\Delta [G(\overline W)]\,=\,0 \quad \textrm{in}\;\;\re^N\setminus\{0\}\,.
\end{equation}
In view of \eqref{e95} there holds
\begin{equation}\label{e93b}
\overline W(x)\ge K \quad \forall x\in \pa B_R \, .
\end{equation}
Furthermore, we have
\begin{equation}\label{e93c}
\overline W(x) \ge a_0 \quad \forall  x\in \pa B_j
\end{equation}
and
\begin{equation}\label{e97}
\overline W(x) \ge G^{-1}\big[G(a_0) +\s\big] \quad \forall x\in \overline{N}_{R,j}\,.
\end{equation}
From \eqref{e93}--\eqref{e97} it follows that $\overline W$ is a
supersolution to the problem
\begin{equation}\label{e98}
\left\{
\begin{array}{ll}
\,   \rho\, \pa_t u = \Delta\big[G(u)\big] &\textrm{in}\,\,N_{R,j}\times (0, T) \\&\\
u\, =\, K &\textrm{on}\,\,\pa B_{R}\times (0, T)\\&\\
u\,=\,a_0 &\textrm{on}\,\, \pa B_j\times(0, T)
\\& \\ \, \! u \, = G^{-1}\big[G(a_0)+\s\big] &\textrm{in}\,\, \overline{N}_{R,j}\times
\{0\}\,.
\end{array}
\right.
\end{equation}
On the other hand, $u_j$ is a subsolution to problem \eqref{e98}.
Hence, by comparison principles,
\begin{equation}\label{e99}
u_j\le \overline W\quad \textrm{in}\;\; N_{R,j}\times (0,T)\,.
\end{equation}

Now let us define
\[\underline W(x):=G^{-1}\big[-\underline{M} \Gamma(x) - \s+ G(a_0) \big] \quad \forall x\in \re^N\setminus\{0\}\,,
\]
where
\[
\underline{M}\ge \frac{G(a_0)-G(-K)}{\Gamma(R)} \,.
\]
By arguments similar to those used above we can infer that $\underline W$ is a subsolution to the problem
\begin{equation}\label{e100}
\left\{
\begin{array}{ll}
\,   \rho\, \pa_t u = \Delta\big[G(u)\big] &\textrm{in}\,\,N_{R,j}\times (0, T) \\&\\
u\, =\,- K &\textrm{on}\,\,\pa B_{R}\times (0, T)\\&\\
u\,=\,a_0 &\textrm{on}\,\, \pa B_j\times(0, T)
\\& \\ \,\! u \, = G^{-1}\big[G(a_0)-\s\big] &\textrm{in}\,\, \overline{N}_{R,j}\times
\{0\}\,.
\end{array}
\right.
\end{equation}
On the other side, $u_j$ is a supersolution to problem
\eqref{e100}. Hence, by comparison principles,
\begin{equation}\label{e101}
\underline{W} \le u_j \quad \textrm{in}\;\; N_{R,j}\times(0,T)\,.
\end{equation}
From \eqref{e99} and \eqref{e101} with $j=j_k$, sending
$k\to\infty$, we obtain
\begin{equation}\label{e102}
\underline W \le u \le \overline W\quad \textrm{in}\;\;
 B_R^c \times (0,T)\,.
\end{equation}
Letting $|x|\to \infty$ in \eqref{e102}, from \eqref{e89} we have that for $|x|$ large enough, \emph{independently of} $t\in [0,T]$, there holds
\[G^{-1}\big[G(a_0)-2\s \big]\le u(x,t) \le G^{-1}\big[G(a_0)+2\s\big] \,. \]
The proof is completed by letting $\s\to 0^+$. \hfill $\square$

\bigskip

In order to prove Theorem \ref{t4} we need some intermediate results.
\begin{lemma}\label{lem1}
Let $N\ge 3$. For any $\alpha,R, M>0$ there exists a subsolution $\underline u_0$ to the equation $-\Delta[G(u)]= 0  \ \textrm{in} \ \re^N$ which is bounded, continuous, radial, nondecreasing as a function of $|x|$, satisfies $\lim_{|x|\to+\infty}\underline u_0(x)=\alpha$ and is
equal to $ -M$ in $B_R$.
\end{lemma}
\noindent \it Proof. \rm
Define
\[\widetilde U_0(x):=G(\alpha) - \frac{\b}{|x|} \quad  \forall x \in B_\e^c \, , \]
where $0<\e <\g:= \frac{\b}{G(\alpha) - G(-M) }$. It is easily seen that ($N \ge 3$)
\[ -\Delta \widetilde U_0(x) \le 0 \quad \forall x \in B_\e^c  \;.\]
Then $\widetilde u_0:=G^{-1}(\widetilde U_0)$ is a subsolution to $-\Delta[G(u)]= 0$ in $B_\e^c$. Consider the function
\[\widehat u_0:=\max\big\{\widetilde u_0,  -M   \big\} \quad \textrm{in} \;\; B_\e^c \,. \]
Since $\widetilde u_0$ solves $-\Delta[G(\widetilde u_0)]\le 0$ in $B_\e^c$, from
Kato's inequality we deduce that $\widehat u_0$ is a subsolution to
$-\Delta[G(u)]= 0$ in $\re^N\setminus \overline B_\e$. Now observe that $\widehat u_0=-M$ in $B_\gamma \setminus B_\e $, so that the function
\[\underline u_0:=
\left\{
\begin{array}{ll}
\,   \widehat u_0 &\textrm{in}\,\, B_\e^c
\\& \\
 -M  & \textrm{in\ \ } B_\e
\end{array}
\right.
\]
is a subsolution to $-\Delta[G(u)]= 0$ in $\re^N$. The fact that $\underline u_0$ is bounded, continuous, radial, nondecreasing as a function of $|x|$ and satisfies the limit property at infinity is clear by construction. The constant condition in $B_R$ is achieved
by choosing $\beta=R\,(G(\alpha)-G(-M))$.  \hfill $\square$

\begin{lemma}\label{lem2}
Suppose that, besides the assumptions of Theorem \ref{t4}, there exists a function $\underline u_0$ having the properties stated in Lemma \ref{lem1} and such that, for a suitable $\e>0$ small enough,
\[
u_0(x)\ge \underline u_0(x)\quad \forall x\in \re^N\,,
\]
\begin{equation}\label{e105bis}
\lim_{|x|\to \infty} \underline u_0(|x|)=\min_{t\in[0,T]}a(t)-\varepsilon>0\,.
\end{equation}
Moreover assume that, for the same $\e$ given above,
\begin{equation}\label{e107}
2\,G\left(\min_{t\in[0,T]}a(t)-\varepsilon\right) > G(\|a\|_\infty)\,.
\end{equation}
Then there exists a solution to problem \eqref{e1}
satisfying condition \eqref{e63}.
\end{lemma}
\noindent \it Proof. \rm  First we repeat the proof of Theorem
\ref{t2} up to the construction of the sequence $\{u_j\}$, keeping
the same notation. Note that, as
in Theorem \ref{t3}, when we allow for a degenerate nonlinearity
$G$, in view of hypothesis $(H_0)$ existence of solutions to
problem \eqref{e68} is due to standard results (see, $e.g.$,
\cite{ACP}). Again, by the results of \cite{DiB3}, $u_j\in C\big(\overline B_j\times [0,T]\big)$.

Then notice that, by the assumptions on $\underline u_0$, \eqref{e105bis}, \eqref{e107} and $(H_0)$, we can find $\beta>0$ and $\widetilde R>\widehat R$ such that for all $R\ge \widetilde
R$
\[\b < \overline u_0(R)\,,\]
\begin{equation}\label{e111}
2 G\big(\overline u_0(R)\big) -  G(\b) - G(\|a\|_\infty)>0\,.
\end{equation}

\smallskip

Still from the assumptions on $\underline u_0$ we deduce that it is a subsolution to problem \eqref{e68}. By comparison principles we have
\[
\underline u_0(|x|)\le u_j(x,t)\le K\quad \forall (x,t) \in
B_j\times(0,T)\,,
\]
where $K$ is as in \eqref{e69}. Hence, by the monotonicity of $\underline u_0$,
\begin{equation}\label{e117}
\underline u_0(R)\le u_j(x,t)\le K \quad \forall (x,t)\in
N_{R,j}\times (0,T)\,.
\end{equation}
Let
\begin{equation}\label{e119}
\g:=\min_{[\b,K]} G' \, .
\end{equation}
Given $\s>0$, in view of \eqref{e71} we can fix
$R=R(\s)>\widetilde R$ in \eqref{e72a} so large that in
\eqref{e72} we are allowed to set
\begin{equation}\label{e112}
\d = \frac 2{\g } V(R)\,.
\end{equation}
Note that $\b$ and $\g$ are independent of $R$ and $\d$. Let
$t_0\in [0,T]$. Define
\begin{equation}\label{e113}
\underline{\l}:= \frac{G(a(t_0))-G(\underline u_0(R))}{\d^2},\quad \underline{M}:=
\frac{2\underline{\l}\d}{\g}\,.
\end{equation}
From \eqref{e111}, \eqref{e112} and \eqref{e113} it follows that
\begin{equation}\label{e114}
\underline{M} = \frac{G(a(t_0))-G(\underline u_0(R))}{V(R)}\,,
\end{equation}
\begin{equation}\label{e115}
-\underline{M} V(R)-\s + G(a(t_0))-\underline{\l} \d^2 > G(\b)
\end{equation}
for $\s>0$ small enough.

Now define
\[\underline w(x,t):=G^{-1}\big[-\underline{M} V(x) -\s +G(a(t_0)) - \underline{\l} (t-t_0)^2\big]   \quad \forall(x,t) \in B_{\widehat R}^c \times  (\underline t_\d, \overline t_\d) \,. \]
Since $|x|\mapsto V(|x|)$ is nonincreasing, by \eqref{e115}
\begin{equation}\label{e118}
\underline w(x,t)\ge \b \quad \forall (x,t)\in
N_{R,j}\times (\underline t_\d, \overline t_\d)\,.
\end{equation}
Also, from $(H_0)(ii)$, \eqref{e73}, \eqref{e118}, \eqref{e119} and
\eqref{e113}
\begin{equation}\label{e74c}
\rho(x)\pa_t \underline w-\Delta[G(\underline
w)]=-\rho(x)\frac{2\underline{\l} (t-t_0)}{G'(\underline w)} + \underline{M}\Delta V\le \rho(x)\left(\frac{2 \underline{\l} \d}{\g} - \underline{M} \right) = 0
\end{equation}
$$ \textrm{in} \;\; B_{\widehat R}^c \times(\underline t_\d,
\overline t_\d) \,. $$
By \eqref{e114},
\begin{equation}\label{e75c}
\underline w(x,t)\le \underline u_0(R) \quad \forall (x,t)\in \pa B_R\times (\underline t_\d , \overline t_\d) \, .
\end{equation}
Furthermore,
\begin{equation}\label{e76c}
\underline w(x,t)\le G^{-1}\big[G(a(t_0)) -\s\big]\quad
\forall (x,t) \in \pa B_j\times (\underline t_\d,
\overline t_\d).
\end{equation}
When $\underline t_\d=0$ there holds
\begin{equation}\label{e78b}
\underline w(x,t)\le G^{-1}\big[G(a(t_0)) -\s\big]\quad \forall (x,t)\in \overline{N}_{R,j}\times \{\underline t_\d\} \, ,
\end{equation}
whereas when $\underline t_\d>0$ we have
\begin{equation}\label{e78c}
\underline w(x,t)\le G^{-1}\big[G(a(t_0)) -\l \d^2\big] =
\underline u_0(R) \quad \forall(x,t) \in \overline{N}_{R,j}\times
\{\underline t_\d\}\,;
\end{equation}
here \eqref{e113} has been used.

From \eqref{e74c}--\eqref{e78c} we infer that $\underline w$ is a subsolution to the problem
\begin{equation}\label{e79c}
\left\{
\begin{array}{ll}
\,   \rho\, \pa_t u = \Delta\big[G(u)\big] &\textrm{in}\,\,N_{R,j}\times (\underline t_\d, \overline t_\d) \\&\\
u\, =\, \underline u_0(R) &\textrm{on}\,\,\pa B_{R}\times (\underline t_\d, \overline t_\d)\\&\\
u\,=\,G^{-1}\big[G(a(t_0))-\s\big] &\textrm{on}\,\, \pa
B_j\times(\underline t_\d, \overline t_\d)
\\& \\ \,\! u \, = \underline u_0(R) &\textrm{in}\,\, \overline{N}_{R,j}\times \{\underline t_\d\}
\end{array}
\right.
\end{equation}
when $\underline t_\d>0$, whereas it is a subsolution to the problem
\begin{equation}\label{e80c}
\left\{
\begin{array}{ll}
\,   \rho\, \pa_t u = \Delta\big[G(u)\big] &\textrm{in}\,\,N_{R,j}\times (\underline t_\d, \overline t_\d) \\&\\
u\, =\, \underline u_0(R) &\textrm{on}\,\,\pa B_{R}\times (\underline t_\d, \overline t_\d)\\&\\
u\,=\,G^{-1}\big[G(a(t_0))-\s\big] &\textrm{on}\,\, \pa
B_j\times(\underline t_\d, \overline t_\d)
\\& \\ \,\! u \, = G^{-1}\big[G(a(t_0))-\s\big] &\textrm{in}\,\, \overline{N}_{R,j}\times \{\underline t_\d\}
\end{array}
\right.
\end{equation}
when $\underline t_\d=0$.

\smallskip
On the other hand, from \eqref{e72}, \eqref{e72a} (which, recall,
holds true as a consequence of \eqref{e65}) and \eqref{e117} we easily
deduce that $u_j$ is a supersolution to problem \eqref{e79c} when
$\underline t_\d>0$, while it is a supersolution to problem
\eqref{e80c} when $\underline t_\d=0$. Hence, by comparison principles,
\[
\underline w\le u_j\quad \textrm{in}\;\; N_{R,j}\times (\underline
t_\d, \overline t_\d)\,.
\]

Finally, let us define
\[\overline w(x,t):=G^{-1}\big[\overline{M} V(x) +\s +G(a(t_0)) +  \overline{\l} (t-t_0)^2\big] \quad \forall(x,t)  \in B_{\widehat R}^c \times (\underline t_\d, \overline t_\d) \, . \]
By construction,
\[\overline w \ge \min_{t \in [0,T]} a(t)  \quad \textrm{in}\;\;  B_{\widehat R}^c \times (\underline t_\d, \overline t_\d) \,.  \]
Choose \[
\overline M\ge \max \left\{\frac{2 \overline{\l} \d}{\min_{t \in [0,T]} G^\prime
(a(t))},
\frac{G(K)-G(-\|a\|_\infty)}{V(R)}\right\},
\]
and
\[
\overline \l \ge \frac{G(K)-G(-\|a\|_\infty)}{\d^2}\,.
\]
Thanks to arguments analogous to those used above, we can infer that
$\overline w$ is a supersolution to problem \eqref{e83} when
$\underline t_\d>0$, whereas it is a supersolution to problem
\eqref{e84} when $\underline t_\d=0$. On the other hand, from \eqref{e72}, \eqref{e72a} and \eqref{e69} we easily deduce that $u_j$ is a subsolution to problem \eqref{e83}
when $\underline t_\d>0$, while it is a subsolution to problem
\eqref{e84} when $\underline t_\d=0$. As in the proof of Theorem \ref{t3}, by means of a compactness argument which makes use of \cite[Lemma 5.2]{DiB1} and a diagonal procedure we deduce that there exists a subsequence $\{u_{j_k}\}\subseteq\{u_j\}$ which
converges, as $k\to \infty,$ locally uniformly in $\re^N\times
(0,T)$ to a solution $u$ to problem \eqref{e1}. We then conclude arguing as in the final part of
the proof of Theorem \ref{t2}. \hfill$\square$

\medskip
\noindent{\it Proof of Theorem \ref{t4}\,.} First consider a
datum $a(t)$ at infinity such that, for some $\e > 0$,
\eqref{e107} holds and $\min_{t \in [0,T]}a(t) - \e > 0$. Consider
then the function $\underline{u}_0$ given in Lemma \ref{lem1} with
the choices $\alpha=\min_{t \in [0,T]}a(t) - \e $, $R$ great
enough so that $u_0(x)\ge \min_{t \in [0,T]}a(t) - \e $ for all $x
\in B_R^c$ and $M = \max(0,-\inf_{x\in \re^N} u_0(x)) $. Clearly,
under these assumptions, $u_0(x)\ge \underline u_0(x)$ for all $x
\in \re^N$. Therefore the assertion of Lemma \ref{lem2} holds
true and the theorem is proved for such $a(\cdot)$.

If there exists no $\e>0$ such that $a(t)$ fulfils \eqref{e107} in the time interval $[0,T]$, we can always find $\e,\tau>0$ small enough such that
\[
2\,G\left(\min_{s\in[t,(t+\tau)\wedge T ]}a(s)-\varepsilon\right) > G\left(\max_{s\in[t,(t+\tau)\wedge T ]}a(s) \right) \ \ \ \forall t \in[0,T)  \,.
\]
This is a consequence of the uniform continuity of $G(a(t))$ and
of its strict positivity in $[0,T]$. Hence we get existence in the
time interval $[0,\tau]$. Repeating this procedure starting from
$t=\tau$ we get existence in the time interval $[\tau,2\tau\wedge
T]$ with initial datum $u(\tau)$ and hence, by Definition
\ref{defsol}, existence in the time interval $[0,2\tau\wedge T]$.
A finite number of iterations yields the claim. \hfill $\square$

\bigskip
\noindent{\it On Remark \ref{texib}.\,} Note that $(i)$ follows from the same proof of
Theorem \ref{t2}, taking $t_0>0$, and $0<\delta <t_0$ in
\eqref{e72}. As for $(ii)$, it is enough to observe that \eqref{er4} permits to repeat the proof of Theorem \ref{t4}, up to choosing $R>R_0$ and $\tau \le \epsilon$.

\section{Uniqueness: proofs}\label{secpuni}
\setcounter{equation}{0}
Let $u_1, u_2$ be any two solutions to problem \eqref{e1}. Define
$$q(x,t):=\left\{
\begin{array}{ll}
 \,   \, \frac{G(u_1)-G(u_2)}{u_1-u_2} &\textrm{if}\,\,u_1(x,t)\neq u_2(x,t)
\\&\\
\textrm{ }0  & \textrm{if}\,\,  u_1(x,t)=u_2(x,t)
\end{array}
\right.  $$
for all $(x,t) \in S_T$. Observe that, in view of $(H_0)(ii)$, $q\ge 0$ in
$S_T$ and $q\in L^\infty(S_T)$. Fix $\t \in (0,T)$. Consider a sequence
$\{q_n\}\subseteq C^\infty(S_T)$ such that for every $n\in \ene$
there hold
\[
\frac 1{n^2} \le q_n\le \|q\|_{L^\infty(S_T)} + \frac
1{n^2}\quad\textrm{in}\;\;\, Q_{n,\t}:=B_n\times(0,\t)\,
\]
and
\begin{equation}\label{e49}
\left\|\frac{(q_n - q)}{\sqrt{q_n}}\right\|_{L^2(Q_{n,\t})} \to 0
\quad \textrm{as}\;\, n\to\infty\,.
\end{equation}
For any $n\in \ene$, let $\psi_n\in C^2(\overline Q_{n,\t})$ be the unique solution to the backward parabolic problem
\begin{equation}\label{e50}
\left\{
\begin{array}{ll}
 \,   \,\rho\, \pa_t \psi_n + q_n\Delta \psi_n = 0 \, &\textrm{in}\,\,Q_{n,\t}
\\&\\
\textrm{ }\psi_n \, = 0& \textrm{on}\,\, \pa B_{n} \times (0,\tau) \\&\\
\textrm{ }\psi_n \, = \chi  &\textrm{in}\,\, \overline{B}_n\times \{\t\}\,,
\end{array}
\right.
\end{equation}
where $\chi\in C^\infty(\re^N),\, 0\le \chi\le 1$ and $supp\;
\chi\subseteq B_{n_0}$ for some fixed $n_0\in \ene$.

\smallskip

The following lemma will play a central role in the proof of
Theorem \ref{t1}.

\begin{lemma}
For every $n\in \ene$ let $\psi_n\in C^2(\overline Q_{n,\t})$ be the
unique solution to problem \eqref{e50}. Then
\begin{equation}\label{e51}
0\le \psi_n \le 1\quad\textnormal{in}\;\;\,Q_{n,\t}\,.
\end{equation}
Furthermore, there exists a constant $C>0$ such that for every
$n>n_0$
\begin{equation}\label{e1b}
-\frac C{n^{N-1}} \, \le \langle \nabla \psi_n, \nu_n \rangle\,
\le 0\quad \textnormal{on}\;\;\pa B_n\times(0,\t),
\end{equation}
where $\nu_n=\nu_n(\sigma)$ is the outer normal at $\sigma \in \pa B_n$\,.
\end{lemma}
\noindent {\it Proof\,.} First notice that $\underline \psi\equiv 0$ is a
subsolution, while $\overline \psi  \equiv 1 $ is a supersolution to problem
\eqref{e50},  so that by comparison we get \eqref{e51}. Now,
since
$$\psi_n=0\quad\textrm{in}\;\;\pa B_n\times (0,\t) $$
for all $n\in \ene$, from \eqref{e51} we deduce that
\begin{equation}\label{e3b}
\langle \nabla \psi_n, \nu_n \rangle\, \le 0\quad
\textrm{in}\;\;\pa B_n\times(0,T)\,.
\end{equation}
For every $n>n_0$ set
\[E_n:= B_n \setminus B_{n_0}\,. \]
From \eqref{e51} and the fact that $\suppo \chi\subset B_{n_0}$ we infer that, for all $n>n_0$, the function $\psi_n$ is a \emph{subsolution} to problem
\begin{equation}\label{e2b}
\left\{
\begin{array}{ll}
 \,   \,\rho\, \pa_t \psi_n + q_n\Delta \psi_n = 0 \,
 &\textrm{in}\,\, E_n \times (0,\tau)
\\&\\
\textrm{ }\psi \, = 1 & \textrm{on}\,\, \pa B_{n_0}\times (0,\tau)  \\&\\
\textrm{ }\psi \, = 0 & \textrm{on}\,\, \pa B_{n}\times (0,\tau) \\&\\
\textrm{ }\psi \, = 0 &\textrm{in}\,\, \overline{B}_n\times \{\t\} \,.
\end{array}
\right.
\end{equation}
Define
\[z(x):= \widehat{C} \, \frac{|x|^{2-N}-n^{2-N}}{1- n^{2-N}}\quad \forall x\in E_n  \, , \]
where $\widehat{C}$ is a positive constant to be chosen. It is easily seen that, for $\widehat C=\widehat C(n_0)$ sufficiently large, the function $z$ is a \emph{supersolution} to problem \eqref{e2b}. Furthermore,
\[\psi_n\,=\,z\,=\,0\quad\textrm{on}\;\; \pa B_n\times(0,\t);\]
hence,
\begin{equation}\label{e4b}
\langle \nabla \psi_n , \nu_n \rangle \ge \langle \nabla  z, \nu_n\rangle= \frac{
(2-N) \, \widehat{C} \, n^{1-N} } {1-n^{2-N}}  \quad \textrm{on}\;\; \pa B_n\times (0,\t)
\end{equation}
for all $n>n_0$. From \eqref{e3b} and \eqref{e4b}, \eqref{e1b} follows with $C:=(N-2)\widehat{C}/(1-n_0^{2-N}).$ This completes the proof.
\hfill $\square$

\medskip

\noindent{\it Proof of Theorem \ref{t1}\,.} Let $u_1$, $u_2$ be two bounded solutions to problem \eqref{e1} satisfying
\[
\lim_{|x|\to \infty} u_1(x,t)=\lim_{|x|\to \infty} u_2(x,t)=a(t) \ \ \ \textnormal{for almost every} \ t \in (0,T) \,.
\]
Clearly, by dominated convergence, this implies that for any $\t\in (0,T)$
\begin{equation}\label{e5b}
\lim_{R\to \infty}\frac 1{R^{N-1}}\int_0^\t\int_{\pa B_R}
\big|G\big(u_1(x,t)\big) -G\big(u_2(x,t)\big)
 \big| \, d\s dt=0\,.
\end{equation}
Put $w:=u_1-u_2\,.$ By Definition \ref{defsol},
\begin{equation}\label{e47}
\begin{split}
\int_{\Omega_1} \rho w(x,\t) \psi(x,\t)dx  = &  \int_0^\t\int_{\Omega_1}\big\{\rho w \psi_t +
[G(u_1)-G(u_2)]\Delta \psi \big\} dx dt  \\
& - \int_0^\t\int_{\pa
\Omega_1} [G(u_1)-G(u_2)]\langle \nabla\psi, \nu\rangle d\s
dt\hspace{.5 cm}
\end{split}
\end{equation}
for any $\t$, $\Omega_1$ and $\psi$ as in Definition \ref{defsol}.

Moreover, multiplying the first equation in \eqref{e50} by ${\Delta\psi_n}/{\rho}$ and integrating by parts we obtain (recall that $\rho\in L^\infty(\re^N)$), for any $n\in \ene$,
\begin{equation}\label{e51bis}
\int_0^\t \int_{B_n} q_n (\Delta\psi_n)^2 dx dt \le \widetilde C
\end{equation}
for some constant $\widetilde C>0$  independent of $n$.

Taking $\Omega_1=B_n$ and $\psi=\psi_n$ in \eqref{e47} we get, for
any $n\in \ene$,
\begin{equation}\label{e52}
\int_{B_n} \rho w(x,\t) \chi dx =
\int_0^\t\int_{B_n}(q-q_n) w \Delta \psi_ndx dt -
\int_0^\t\int_{\pa B_n} q w \langle\nabla \psi_n , \nu_n\rangle \, d\s dt.
\end{equation}
We shall prove that both integrals on the right-hand side of inequality \eqref{e52} tend to $0$ as $n\to \infty$. In fact, from \eqref{e49} and \eqref{e51bis} we have:
\begin{equation}\label{e53}
\begin{split}
& \left(\int_0^\t\int_{B_n}(q-q_n) w \Delta\psi_ndx
dt\right)^2  \\
\le &  \overline{C} \int_0^\t\int_{B_n}\left|\frac{q-q_n}{\sqrt{q_n}}\right|^2dx
dt\int_0^\t \int_{B_n}q_n |\Delta \psi_n|^2 dx dt \to
0\quad \textrm{as} \;\; n\to \infty \, ,
\end{split}
\end{equation}
where $\overline C:=(\|u_1\|_\infty + \|u_2\|_\infty)^2$.

Moreover, by \eqref{e1b} and \eqref{e5b}, for every $n>n_0$ there holds
\begin{equation}\label{e6b}
\begin{split}
 \left| \int_0^\t \int_{\pa B_n} q w \langle\nabla \psi_n, \nu_n
\rangle d\s dt \right| = &  \left|\int_0^\t \int_{\pa
B_n}\big[G(u_1)-G(u_2)\big] \langle\nabla \psi_n, \nu_n \rangle
d\s dt \right|  \\
\le & \max_{\pa B_n}|\langle\nabla \psi_n, \nu_n \rangle|\int_0^\t
\int_{\pa B_n}\big|G(u_1)-G(u_2)\big|d\s dt \\
\le & \frac{C}{n^{N-1}}\int_0^\t \int_{\pa B_n}
\big|G(u_1)-G(u_2)\big|d\s dt\to 0
\end{split}
\end{equation}
as $n\to \infty$. Sending $n\to \infty$ in \eqref{e52}, from \eqref{e53} and
\eqref{e6b} it follows that
\begin{equation}\label{e60}
\int_{\re^N} \rho(x) \chi(x) w(x,\t)dx = 0
\end{equation}
for any $\t\in(0,T)$ and any $\chi\in C^\infty_c(\re^N)$ with $0\le \chi\le 1$.

\smallskip
Now fix any compact subset $K\subset \re^N$.  Define
\[\zeta(x,\tau):=\left\{
\begin{array}{ll}
 \,   \, 1 &\textrm{if} \,\, x\in K \,\, \textnormal{and}  \,\, w(x,\tau)>0 \, ,
\\&\\
\textrm{ }0  & \textrm{elsewhere} \, .
\end{array}
\right.
\]
Pick a sequence $\{\chi_n\}\subseteq C^\infty_c(\re^N)$, with $0\le \chi_n\le 1$, such that
$\chi_n(x)\to \zeta(x)$ as $n\to\infty$ for any $x\in \re^N$. In view of \eqref{e60} we deduce that, for any $n\in \ene$,
\begin{equation}\label{e61}
\int_{\re^N } \rho(x) \chi_n(x) w(x,\t) \, dx = 0 \, .
\end{equation}
Letting $n\to \infty$ in \eqref{e61}, by dominated convergence we get
\[
\int_{K  \cap \{ w(\cdot,\tau) > 0 \} } \rho(x)w(x,\t) dx=0\,.
\]
Hence $w(x,\t) \le  0$ for any $x\in K$. Since the compact subset $K\subset \re^N$ and $\t\in (0,T)$ are arbitrary, we get
\[w \le  0\quad \textrm{in}\;\; \re^N\times (0,T) \, , \]
that is
\[u_1 \le u_2\quad \textrm{in}\;\; \re^N\times (0,T) \, .\]
Interchanging the role of $u_1$ and $u_2$ we obtain also the opposite inequality, and this completes the proof. \hfill $\square$

\end{document}